\documentclass[12pt,draftcls,onecolumn]{IEEEtran}

\usepackage{cite}
\usepackage{graphicx}
\usepackage[breaklinks=false, allcolors=black, colorlinks=true]{hyperref}

\usepackage{amsfonts}
\usepackage{amsmath} 
\usepackage{amssymb}
\usepackage{amsthm}	
\usepackage{enumerate}
\usepackage{mathrsfs}
\usepackage{subcaption}
\usepackage{epstopdf}
\usepackage{dashrule}
\usepackage[normalem]{ulem}

\theoremstyle{definition} 
\newtheorem{theore}{Theorem}
\newtheorem{corollar}{Corollary}
\newtheorem{lemm}{Lemma}
\newtheorem{definitio}{Definition}
\newtheorem{proble}{Problem}
\newtheorem{assumptio}{Assumption}
\newtheorem{propert}{Property}
\newtheorem{propositio}{Proposition}
\newtheorem{conjectur}{Conjecture}
\newtheorem{exampl}{Example}
\newtheorem{rul}{Rule}

\newtheorem{remar}{Remark}
\newtheorem{fac}{Fact}

\def\nn{\nonumber}
\def\beq{\begin{equation}}
\def\eeq{\end{equation}}
\def\beqs{\begin{equation*}}
\def\eeqs{\end{equation*}}
\def\balign{\begin{align}}
\def\ealign{\end{align}}
\def\bea{\begin{eqnarray}}
\def\eea{\end{eqnarray}}
\def\ba{\begin{array}}
\def\ea{\end{array}}
\def\bcenter{\begin{center}}
\def\ecenter{\end{center}}
\def\bitem{\begin{itemize}}
\def\eitem{\end{itemize}}
\def\benum{\begin{enumerate}}
\def\eenum{\end{enumerate}}
\def\bdoc{\begin{document}}
\def\edoc{\end{document}}
\def\defeq{{\stackrel{\Delta}{=}}}
\def\nn{\nonumber}

\def\eg{{e.g., \/}}
\def\etal{{\it et al. \/}}
\def\viz{{\iet viz.,\ \/}}
\def\ie{{i.e.,\ \/}}
\def\etc{{\it etc.,\ \/}}
\def\cf{{cf.,\ \/}}

\usepackage[usenames,dvipsnames,svgnames,table]{xcolor}

\def\yellow{\textcolor{yellow}}
\def\green{\textcolor{green}}
\def\tcbg{\textcolor{bgrd}}
\def\magenta{\textcolor{magenta}}
\def\white{\textcolor{white}}
\def\gray{\textcolor{gray}}
\def\cyan{\textcolor{cyan}}

\newcommand\Quote[1]{\lq\textsl{#1}\rq}
\newcommand\fr[2]{{\textstyle\frac{#1}{#2}}}

\newcommand{\Ambb}{\mathbb{A}}
\newcommand{\Bmbb}{\mathbb{B}}
\newcommand{\Cmbb}{\mathbb{C}}
\newcommand{\Dmbb}{\mathbb{D}}
\newcommand{\Embb}{\mathbb{E}}
\newcommand{\Fmbb}{\mathbb{F}}
\newcommand{\Gmbb}{\mathbb{G}}
\newcommand{\Hmbb}{\mathbb{H}}
\newcommand{\Imbb}{\mathbb{I}}
\newcommand{\Jmbb}{\mathbb{J}}
\newcommand{\Kmbb}{\mathbb{K}}
\newcommand{\Lmbb}{\mathbb{L}}
\newcommand{\Mmbb}{\mathbb{M}}
\newcommand{\Nmbb}{\mathbb{N}}
\newcommand{\Ombb}{\mathbb{O}}
\newcommand{\Pmbb}{\mathbb{P}}
\newcommand{\Qmbb}{\mathbb{Q}}
\newcommand{\Rmbb}{\mathbb{R}}
\newcommand{\Smbb}{\mathbb{S}}
\newcommand{\Tmbb}{\mathbb{T}}
\newcommand{\Umbb}{\mathbb{U}}
\newcommand{\Vmbb}{\mathbb{V}}
\newcommand{\Wmbb}{\mathbb{W}}
\newcommand{\Xmbb}{\mathbb{X}}
\newcommand{\Ymbb}{\mathbb{Y}}
\newcommand{\Zmbb}{\mathbb{Z}}

\newcommand{\Amsc}{\mathcal{A}}
\newcommand{\Bmsc}{\mathcal{B}}
\newcommand{\Cmsc}{\mathcal{C}}
\newcommand{\Dmsc}{\mathcal{D}}
\newcommand{\Emsc}{\mathcal{E}}
\newcommand{\Fmsc}{\mathcal{F}}
\newcommand{\Gmsc}{\mathcal{G}}
\newcommand{\Hmsc}{\mathcal{H}}
\newcommand{\Imsc}{\mathcal{I}}
\newcommand{\Jmsc}{\mathcal{J}}
\newcommand{\Kmsc}{\mathcal{K}}
\newcommand{\Lmsc}{\mathcal{L}}
\newcommand{\Mmsc}{\mathcal{M}}
\newcommand{\Nmsc}{\mathcal{N}}
\newcommand{\Omsc}{\mathcal{O}}
\newcommand{\Pmsc}{\mathcal{P}}
\newcommand{\Qmsc}{\mathcal{Q}}
\newcommand{\Rmsc}{\mathcal{R}}
\newcommand{\Smsc}{\mathcal{S}}
\newcommand{\Tmsc}{\mathcal{T}}
\newcommand{\Umsc}{\mathcal{U}}
\newcommand{\Vmsc}{\mathcal{V}}
\newcommand{\Wmsc}{\mathcal{W}}
\newcommand{\Xmsc}{\mathcal{X}}
\newcommand{\Ymsc}{\mathcal{Y}}
\newcommand{\Zmsc}{\mathcal{Z}}

\newcommand{\Ascr}{\mathscr{A}}
\newcommand{\Bscr}{\mathscr{B}}
\newcommand{\Cscr}{\mathscr{C}}
\newcommand{\Dscr}{\mathscr{D}}
\newcommand{\Escr}{\mathscr{E}}
\newcommand{\Fscr}{\mathscr{F}}
\newcommand{\Gscr}{\mathscr{G}}
\newcommand{\Hscr}{\mathscr{H}}
\newcommand{\Iscr}{\mathscr{I}}
\newcommand{\Jscr}{\mathscr{J}}
\newcommand{\Kscr}{\mathscr{K}}
\newcommand{\Lscr}{\mathscr{L}}
\newcommand{\Mscr}{\mathscr{M}}
\newcommand{\Nscr}{\mathscr{N}}
\newcommand{\Oscr}{\mathscr{O}}
\newcommand{\Pscr}{\mathscr{P}}
\newcommand{\Qscr}{\mathscr{Q}}
\newcommand{\Rscr}{\mathscr{R}}
\newcommand{\Sscr}{\mathscr{S}}
\newcommand{\Tscr}{\mathscr{T}}
\newcommand{\Uscr}{\mathscr{U}}
\newcommand{\Vscr}{\mathscr{V}}
\newcommand{\Wscr}{\mathscr{W}}
\newcommand{\Xscr}{\mathscr{X}}
\newcommand{\Yscr}{\mathscr{Y}}
\newcommand{\Zscr}{\mathscr{Z}}

\def\Atiny{\mbox{\tiny{A}}}
\def\Btiny{\mbox{\tiny{B}}}
\def\Ctiny{\mbox{\tiny{C}}}
\def\Dtiny{\mbox{\tiny{D}}}
\def\Etiny{\mbox{\tiny{E}}}
\def\Ctiny{\mbox{\tiny{C}}}
\def\Ftiny{\mbox{\tiny{F}}}
\def\Gtiny{\mbox{\tiny{G}}}
\def\Htiny{\mbox{\tiny{H}}}
\def\Itiny{\mbox{\tiny{I}}}
\def\Jtiny{\mbox{\tiny{J}}}
\def\Ktiny{\mbox{\tiny{K}}}
\def\Ltiny{\mbox{\tiny{L}}}
\def\Mtiny{\mbox{\tiny{M}}}
\def\Ntiny{\mbox{\tiny{N}}}
\def\Otiny{\mbox{\tiny{O}}}
\def\Ptiny{\mbox{\tiny{P}}}
\def\Qtiny{\mbox{\tiny{Q}}}
\def\Rtiny{\mbox{\tiny{R}}}
\def\Stiny{\mbox{\tiny{S}}}
\def\Ttiny{\mbox{\tiny{T}}}
\def\Utiny{\mbox{\tiny{U}}}
\def\Vtiny{\mbox{\tiny{V}}}
\def\Wtiny{\mbox{\tiny{W}}}
\def\Xtiny{\mbox{\tiny{X}}}
\def\Ytiny{\mbox{\tiny{Y}}}
\def\Ztiny{\mbox{\tiny{Z}}}

\def\alphabf{\hbox{\boldmath$\alpha$\unboldmath}}
\def\betabf{\hbox{\boldmath$\beta$\unboldmath}}
\def\gammabf{\hbox{\boldmath$\gamma$\unboldmath}}
\def\deltabf{\hbox{\boldmath$\delta$\unboldmath}}
\def\epsilonbf{\hbox{\boldmath$\epsilon$\unboldmath}}
\def\zetabf{\hbox{\boldmath$\zeta$\unboldmath}}
\def\etabf{\hbox{\boldmath$\eta$\unboldmath}}
\def\iotabf{\hbox{\boldmath$\iota$\unboldmath}}
\def\kappabf{\hbox{\boldmath$\kappa$\unboldmath}}
\def\lambdabf{\hbox{\boldmath$\lambda$\unboldmath}}
\def\mubf{\hbox{\boldmath$\mu$\unboldmath}}
\def\nubf{\hbox{\boldmath$\nu$\unboldmath}}
\def\xibf{\hbox{\boldmath$\xi$\unboldmath}}
\def\pibf{\hbox{\boldmath$\pi$\unboldmath}}
\def\rhobf{\hbox{\boldmath$\rho$\unboldmath}}
\def\sigmabf{\hbox{\boldmath$\sigma$\unboldmath}}
\def\taubf{\hbox{\boldmath$\tau$\unboldmath}}
\def\upsilonbf{\hbox{\boldmath$\upsilon$\unboldmath}}
\def\phibf{\hbox{\boldmath$\phi$\unboldmath}}
\def\chibf{\hbox{\boldmath$\chi$\unboldmath}}
\def\psibf{\hbox{\boldmath$\psi$\unboldmath}}
\def\omegabf{\hbox{\boldmath$\omega$\unboldmath}}
\def\Sigmabf{\hbox{$\bf \Sigma$}}
\def\Upsilonbf{\hbox{$\bf \Upsilon$}}
\def\Omegabf{\hbox{$\bf \Omega$}}
\def\Deltabf{\hbox{$\bf \Delta$}}
\def\Gammabf{\hbox{$\bf \Gamma$}}
\def\Thetabf{\hbox{$\bf \Theta$}}
\def\Lambdabf{\mbox{$ \bf \Lambda $}}
\def\Xibf{\hbox{\bf$\Xi$}}
\def\Pibf{{\bf \Pi}}
\def\thetabf{{\mbox{\boldmath$\theta$\unboldmath}}}
\def\Upsilonbf{\hbox{\boldmath$\Upsilon$\unboldmath}}
\newcommand{\Phibf}{\mbox{${\bf \Phi}$}}
\newcommand{\Psibf}{\mbox{${\bf \Psi}$}}

\def\abf{{\bf a}}
\def\bbf{{\bf b}}
\def\cbf{{\bf c}}
\def\dbf{{\bf d}}
\def\ebf{{\bf e}}
\def\fbf{{\bf f}}
\def\gbf{{\bf g}}
\def\hbf{{\bf h}}
\def\ibf{{\bf i}}
\def\jbf{{\bf j}}
\def\kbf{{\bf k}}
\def\lbf{{\bf l}}
\def\mbf{{\bf m}}
\def\nbf{{\bf n}}
\def\obf{{\bf o}}
\def\pbf{{\bf p}}
\def\qbf{{\bf q}}
\def\rbf{{\bf r}}
\def\sbf{{\bf s}}
\def\tbf{{\bf t}}
\def\ubf{{\bf u}}
\def\vbf{{\bf v}}
\def\wbf{{\bf w}}
\def\xbf{{\bf x}}
\def\ybf{{\bf y}}
\def\zbf{{\bf z}}
\def\rbf{{\bf r}}
\def\xbf{{\bf x}}
\def\ybf{{\bf y}}
\def\Abf{{\bf A}}
\def\Bbf{{\bf B}}
\def\Cbf{{\bf C}}
\def\Dbf{{\bf D}}
\def\Ebf{{\bf E}}
\def\Fbf{{\bf F}}
\def\Gbf{{\bf G}}
\def\Hbf{{\bf H}}
\def\Ibf{{\bf I}}
\def\Jbf{{\bf J}}
\def\Kbf{{\bf K}}
\def\Lbf{{\bf L}}
\def\Mbf{{\bf M}}
\def\Nbf{{\bf N}}
\def\Obf{{\bf O}}
\def\Pbf{{\bf P}}
\def\Qbf{{\bf Q}}
\def\Rbf{{\bf R}}
\def\Sbf{{\bf S}}
\def\Tbf{{\bf T}}
\def\Ubf{{\bf U}}
\def\Vbf{{\bf V}}
\def\Wbf{{\bf W}}
\def\Xbf{{\bf X}}
\def\Ybf{{\bf Y}}
\def\Zbf{{\bf Z}}
\def\Ac{{\cal A}}
\def\Bc{{\cal B}}
\def\Cc{{\cal C}}
\def\Dc{{\cal D}}
\def\Ec{{\cal E}}
\def\Fc{{\cal F}}
\def\Gc{{\cal G}}
\def\Hc{{\cal H}}
\def\Ic{{\cal I}}
\def\Jc{{\cal J}}
\def\Kc{{\cal K}}
\def\Lc{{\cal L}}
\def\Mc{{\cal M}}
\def\Nc{{\cal N}}
\def\Oc{{\cal O}}
\def\Pc{{\cal P}}
\def\Qc{{\cal Q}}
\def\Rc{{\cal R}}
\def\Sc{{\cal S}}
\def\Tc{{\cal T}}
\def\Uc{{\cal U}}
\def\Vc{{\cal V}}
\def\Wc{{\cal W}}
\def\Xc{{\cal X}}
\def\Yc{{\cal Y}}
\def\Zc{{\cal Z}}

\def\0bf{\mathbf{0}}
\def\1bf{\mathbf{1}}

\def\lambdabf{\boldsymbol{\lambda}}
\def\mubf{\boldsymbol{\mu}}
\def\gammabf{\boldsymbol{\gamma}}
\def\xibf{\boldsymbol{\xi}}
\def\sigmabf{\boldsymbol{\sigma}}
\def\zetabf{\boldsymbol{\zeta}}

\def\Lambdabf{\boldsymbol{\Lambda}}
\def\Mubf{\boldsymbol{\Mu}}
\def\Gammabf{\boldsymbol{\Gamma}}
\def\Xibf{\boldsymbol{\Xi}}
\def\Sigmabf{\boldsymbol{\Sigma}}
\def\Zetabf{\boldsymbol{\Zeta}}

\newcommand\independent{\protect\mathpalette{\protect\independenT}{\perp}}
\def\independenT#1#2{\mathrel{\rlap{$#1#2$}\mkern2mu{#1#2}}}

\newcommand{\bsf}[1]{\textsf{\textbf{#1}}}
\newcommand{\lbsf}[1]{\textsf{\large  \textbf{#1}}}
\newcommand{\Lbsf}[1]{\textsf{\Large  \textbf{#1}}}
\newcommand{\hbsf}[1]{\textsf{\huge  \textbf{#1}}}

\newcommand{\myminipage}[3]{\begin{minipage}[#1]{#2}{#3} \end{minipage}}
\newcommand{\sbs}[4]{\myminipage{c}{#1}{#3} \hfill
\myminipage{c}{#2}{#4}}

\newcommand{\myfig}[2]{\centerline{\psfig{figure=#1,width=#2,silent=}}}
\newcommand{\myfigh}[2]{\centerline{\psfig{figure=#1,height=#2,silent=}}}
\newcommand{\myfigwh}[3]{\centerline{\psfig{figure=#1,width=#2,height=#3,silent=}}}

\newcommand{\beqa}{\begin{eqnarray}}
\newcommand{\eeqa}{\end{eqnarray}}
\newcommand{\beqan}{\begin{eqnarray*}}
\newcommand{\eeqan}{\end{eqnarray*}}
\newcommand{\dst}[1]{\displaystyle{ #1 }}
\newcounter{pf}

\newcommand{\smax}[1] { \bar \sigma \left( #1 \right) }
\newcommand{\Rn}{{\mathbb R}^n}
\newcommand{\R}{{\mathbb R}}
\newcommand{\C}{{\mathbb C}}
\newcommand{\Rm}{\mathbb{R}^m}
\newcommand{\Rmn}{\mathbb{R}^{m \times n}}
\newcommand{\Rpq}{\mathbb{R}^{p \times q}}
\newcommand{\Cn}{\mathbb{C}^n}
\newcommand{\Cm}{\mathbb{C}^m}
\newcommand{\Cnn}{\mathbb{C}^{n \times n}}
\newcommand{\Cmn}{\mathbb{C}^{m \times n}}
\newcommand{\ip}[1]{\left\langle #1 \right\rangle}
\newcommand{\rank}{\mbox{rank}}
\newcommand{\Span}{\mbox{\rm Span }}
\newcommand{\Trace}{\mbox{\rm Tr }}
\newcommand{\trace}[1]{\text{Tr}\left(#1\right)}
\newcommand{\Spec}{\mbox{\rm Spec }}
\newcommand{\vectornorm}[1]{\left\|#1\right\|}

\newcommand{\pd}[2]{\frac{\partial #1}{\partial #2}}
\newcommand{\ppd}[3]{\frac{\partial^2 #1}{\partial #2 \partial #3}}

\newcommand{\thtilde}{\tilde{\theta}}
\newcommand{\thnom}{\theta^\circ}
\newcommand{\thopt}{\theta^{\mbox{\small opt}}}
\newcommand{\thhat}{{\hat{\theta}}}
\newcommand{\Tho}{\Theta^\circ}
\newcommand{\tho}{\theta^\circ}
\newcommand{\np}{{n_p}}

\newcommand{\ii}{{[i]}}
\newcommand{\II}{{[i+1]}}
\newcommand{\iii}{{[ii]}}
\newcommand{\jj}{{[j]}}
\newcommand{\kk}{{[k]}}
\newcommand{\thi}{{\theta^\ii}}
\newcommand{\thI}{{\theta^\II}}
\newcommand{\di}{{d^\ii}}
\newcommand{\gi}{{g^\ii}}
\newcommand{\Hi}{{\HH^\ii}}
\newcommand{\thK}{\theta^{(k+1)}}
\newcommand{\gk}{{g^{(k)}}}
\newcommand{\Hk}{{{\cal H}^{(k)}}}

\newcommand{\bfdelta}{{\bf \Delta}}

\newcommand{\Exp}[1]{\exp \left\{ #1 \right\}} 
\newcommand{\gaussian}[1]{\mathbb{N} \left( #1 \right)}
\newcommand{\uniform}[1]{\mathbb{U} \left[ #1 \right]}
\newcommand{\exponential}[1]{\mathbb{E} \left[ #1 \right]}
\newcommand{\EXP}[1]{\EEXP \left[ #1 \right]} 
\newcommand{\EEXP}{\mbox{\bsf{E}}} 
\newcommand{\Prob}[1]{\mbox{{\sf Pr}} \left(#1 \right)}
\newcommand{\convas}{\stackrel{as}{\longrightarrow}}
\newcommand{\convinp}{\stackrel{p}{\longrightarrow}}
\newcommand{\convind}{\stackrel{d}{\longrightarrow}}
\newcommand{\convqm}{\stackrel{qm}{\longrightarrow}}
\newcommand{\sss}[1]{{_{#1}}}
\newcommand{\density}[2]{p_{_{_{#1}}}\!\!\left(#2 \right)} 
\newcommand{\distro}[2]{P_{_{_{#1}}}\!\!\left(#2 \right)} 
\newcommand{\rxx}[1]{R_{_{#1}}\!} 
\newcommand{\sxx}[1]{S_{_{#1}}} 
\newcommand{\cov}[1]{\Lambda_{_{#1}}} 
\newcommand{\mean}[1]{m_{_{#1}}} 
\newcommand{\LS}[1]{\hat{#1}_{_{LS}}} 
\newcommand{\MV}[1]{\hat{#1}_{_{MV}}} 
\newcommand{\LMV}[1]{\hat{#1}_{_{LMV}}} 
\newcommand{\ML}[1]{\hat{#1}_{_{ML}}} 

\renewcommand{\arraystretch}{0.9}
\newcommand{\bmat}[1]{ \begin{bmatrix} #1 \end{bmatrix}}
\newcommand{\mat}[1]{ \left[ \begin{array}{cccccccc} #1 \end{array}
\right] }
\newcommand{\smallmat}[1]{\small{\mat{#1}}}
\newcommand{\sysblk}[4]{\begin{array}{c|cccc}#1&#2\\ \hline#3&#4
\end{array}}
\newcommand{\sysmat}[4]{\left[\sysblk{#1}{#2}{#3}{#4}\right]}
\newcommand{\SGeq}{\succ}
\newcommand{\SLeq}{\prec}
\newcommand{\Geq}{\succeq}
\newcommand{\Leq}{\preceq}

\newcommand{\Bset}{\mathbb{B}}
\newcommand{\Cset}{\mathbb{C}}
\newcommand{\Fset}{\mathbb{F}}
\newcommand{\Mset}{\mathbb{M}}
\newcommand{\Nset}{\mathbb{N}}
\newcommand{\Qset}{\mathbb{Q}}
\newcommand{\Rset}{\mathbb{R}}
\newcommand{\Sset}{\mathbb{S}}
\newcommand{\Tset}{\mathbb{T}}
\newcommand{\Uset}{\mathbb{U}}
\newcommand{\Vset}{\mathbb{V}}
\newcommand{\Wset}{\mathbb{W}}
\newcommand{\Zset}{\mathbb{Z}}

\newcommand{\Ical}{{\cal I}}
\newcommand{\Acal}{{\cal A}}
\newcommand{\Bcal}{{\cal B}}
\newcommand{\Ccal}{{\cal C}}
\newcommand{\Dcal}{{\cal D}}
\newcommand{\Ecal}{{\cal E}}
\newcommand{\Fcal}{{\cal F}}
\newcommand{\Gcal}{{\cal G}}
\newcommand{\Hcal}{{\cal H}}
\newcommand{\Kcal}{{\cal K}}
\newcommand{\Lcal}{{\cal L}}
\newcommand{\Mcal}{{\cal M}}
\newcommand{\Ncal}{{\cal N}}
\newcommand{\Pcal}{{\cal P}}
\newcommand{\Qcal}{{\cal Q}}
\newcommand{\Rcal}{{\cal R}}
\newcommand{\Scal}{{\cal S}}
\newcommand{\Tcal}{{\cal T}}
\newcommand{\Wcal}{{\cal W}}
\newcommand{\Ucal}{{\cal U}}
\newcommand{\Vcal}{{\cal V}}
\newcommand{\Xcal}{{\cal X}}
\newcommand{\Zcal}{{\cal Z}}

\newcommand{\EE}{{\bf E}}
\newcommand{\FF}{{\bf F}}
\newcommand{\GG}{{\bf G}}
\newcommand{\HH}{{\bf H}}
\newcommand{\LL}{{\bf L}}
\newcommand{\NN}{{\bf N}}
\newcommand{\MM}{{\bf M}}
\newcommand{\PP}{{\bf P}}
\newcommand{\QQ}{{\bf Q}}
\newcommand{\RR}{{\bf R}}
\renewcommand{\SS}{{\bf S}}
\newcommand{\TT}{{\bf T}}
\newcommand{\VV}{{\bf V}}
\newcommand{\WW}{{\bf W}}

\newcommand{\thk}{\theta^{(k)}}
\newcommand{\thb}{\theta^{\rm opt}}
\newcommand{\alb}{\alpha^{\rm opt}}
\newcommand{\dk}{d^{(k)}}
\newcommand{\Hinf}{{\cal H}_\infty}
\newcommand{\Htwo}{{\cal H}_2}

\renewcommand{\arraystretch}{1.1}

\newcommand{\red}[1]{{\color{red} #1}}
\newcommand{\blue}[1]{{\color{Blue} #1}}
\newcommand{\black}[1]{\color{black} #1}

\newcounter{l1}
\newcounter{l2}
\newcounter{l3}
\newcommand{\bdotlist}{\begin{list}{$\bullet$}{}}
\newcommand{\bboxlist}{\begin{list}{$\Box$}{}}
\newcommand{\bbboxlist}{\begin{list}{\raisebox{.005in}{{\tiny
$\blacksquare$ \ \ }}}{}}
\newcommand{\bdashlist}{\begin{list}{$-$}{} }
\newcommand{\blist}{\begin{list}{}{} }
\newcommand{\barablist}{\begin{list}{\arabic{l1}}{\usecounter{l1}}}
\newcommand{\balphlist}{\begin{list}{(\alph{l2})}{\usecounter{l2}}}
\newcommand{\bAlphlist}{\begin{list}{\Alph{l2}.}{\usecounter{l2}}}
\newcommand{\bdiamlist}{\begin{list}{$\diamond$}{}}
\newcommand{\bromalist}{\begin{list}{(\roman{l3})}{\usecounter{l3}}}

\newcommand{\thm}[1]{\noindent \begin{theorem} #1   \end{theorem}}
\newcommand{\prop}[1]{\begin{proposition} #1 \end{proposition}}
\newcommand{\lem}[1]{\begin{lemma} #1  \hfill $\blacksquare$ \end{lemma}}
\newcommand{\ex}[1]{\begin{example} {\rm #1} \end{example}}
\newcommand{\prf}[1]{ \noindent {\em Proof:} \, #1 \hfill $\blacksquare$}
\newcommand{\rem}[1]{\begin{remark} {\rm #1} \hfill $\Box$ \end{remark}}
\newcommand{\defn}[1]{\begin{definition} {\rm #1 } \end{definition}}
\newcommand{\prob}[1]{\begin{exercise} {\rm  #1 } \end{exercise}}
\newcommand{\cor}[1]{\begin{corollary}   #1  \end{corollary}}

\usepackage{tikz}
\usetikzlibrary{arrows,snakes,backgrounds}
\usepackage{pgfplots}
\pgfplotsset{compat=newest}
\usetikzlibrary{plotmarks}
\usepackage{grffile}

\newcommand{\argmin}{\mathop{\rm argmin}}
\newcommand{\argmax}{\mathop{\rm argmax}}
\newcommand{\diag}{\mathop{\mathrm{diag}}}
\newcommand{\tr}{\mathop{\rm Tr}}
\newcommand{\conv}{\mathop{\rm conv}}
\newcommand{\var}{\mathop{\rm Var}}
\renewcommand{\b}[1]{\ensuremath{\boldsymbol{\mathrm{#1}}}}
\newcommand{\ms}{{\rm MS}}
\newcommand{\tcs}{{\rm TCS}}
\newcommand{\scs}{{\rm SCS}}
\newcommand{\trans}{^{\sf T}}
\newcommand{\E}[1]{\b{\mu}_{{#1}}}
\newcommand{\Var}[1]{{\Sigma_{#1}}}

\newcommand{\bone}{\mathbf{1}}

\def\u{u}	
\def\U{U}	
\def\v{v}	
\def\V{V}	
\def\s{s}	
\def\S{S}	
\def\f{q}	
\def\q{q}	
\def\Q{Q}	
\def\c{c}	
\def\z{z}	
\def\xcap{b}
\def\ccap{\mathbf{\c}}
\def\rentftr{\Phi}
\def\rentfsr{\Sigma}
\def\rentDer{\Delta}
\def\lambdao{\lambda^{o}}
\def \col{{\rm Col}}

\def \subvec{r}

\def \add [#1]{\blue{#1}}
\def \replace [#1]#2{\red{#1} \blue{#2}}

\newcommand{\remove}[1]{}

\def \wl [#1]{\magenta{(\textbf{Sam says:} #1)}}
\def \eb [#1]{\red{(\textbf{Eilyan says:} #1)}}

\def \rone{\black}
\def \rtwo{\black}
\def \rsam{\black}
\def \rthree{\black}

\graphicspath{{figures/}}

\begin{document}

\title{A Convex Information Relaxation for Constrained Decentralized Control Design Problems}

\author{\vspace{.12in} Weixuan Lin   \qquad Eilyan Bitar 
\thanks{Supported in part by NSF grant ECCS-1351621, NSF grant IIP-1632124, US DOE under the CERTS initiative, and the Simons Institute for the Theory of Computing. This technical note builds on the authors' preliminary results published as part of the 2016 IEEE American Control Conference (ACC) \cite{Lin2016}.}
\thanks{W. Lin (wl476@cornell.edu) and E. Bitar (eyb5@cornell.edu) are with the School of Electrical and Computer Engineering, Cornell University, Ithaca, NY, 14853, USA. }
\vspace{-.26in}
}

\maketitle

\begin{abstract}
We describe a convex programming approach to the calculation of lower bounds on the minimum cost of constrained decentralized control problems with nonclassical information structures. The  class of problems we consider entail the decentralized output feedback control of a linear time-varying system over a finite horizon, subject to polyhedral constraints on the state and input trajectories, and sparsity constraints on  the controller's information structure. \rsam{As the determination of} optimal control policies for such systems is known to be computationally intractable in general, considerable effort has been made in the literature to identify efficiently computable, albeit suboptimal, feasible control policies. The construction of computationally tractable bounds on their suboptimality \rsam{is} the primary motivation for the techniques developed in this note. \rone{Specifically, given a decentralized control problem with nonclassical information, we characterize an expansion of the given information structure, which ensures its partial nestedness, while maximizing the optimal value of the resulting decentralized control problem  under the expanded information structure.} The resulting decentralized control problem is cast as an infinite-dimensional convex program, which is further relaxed  via a partial dualization and restriction to affine dual control policies. The resulting problem is a finite-dimensional conic program whose optimal value is a provable lower bound on the minimum cost of the original constrained decentralized control problem. 
\end{abstract}

\section{Introduction} \label{sec:introduction}

In general, the design of an optimal decentralized controller amounts to an infinite-dimensional nonconvex optimization problem.
The difficulty in solution derives in part from the manner in which information is shared between controllers---the so called information structure of a problem; see \cite{Mahajan2012_survey} for a survey.
Considerable effort has been made to identify information structures under which the problem of decentralized control design can be recast as an equivalent \rsam{convex program.}
For instance, partial nestedness of the information structure \cite{Ho1972} is known to simplify the control design, as it eliminates the incentive to signal between controllers. In particular, linear control policies are guaranteed to be optimal for decentralized LQG problems with partially nested information structures \cite{Ho1972}.
Closely related notions of quadratic invariance \cite{Rotkowitz2006} and funnel causality \cite{Bamieh2005} guarantee convexity of decentralized controller synthesis,  
whose objective is to minimize the closed-loop norm of an LTI system. \remove{There is also a body of literature offering insight on the structure of optimal decentralized controllers for problems with  nonclassical information structures; we refer the reader to   \cite{Mahajan2014_partial_history, Nayyar2011_delayed_sharing, Nayyar2013_partial_history, sandell1974solution} for recent advances.}

\rone{As the tractable computation of optimal policies for the majority of decentralized control problems with nonclassical information structures remains out of reach \cite{Mahajan2012_survey}, there is a practical need to quantify the suboptimality of feasible policies via the derivation of lower bounds on the optimal values of such problems. \rsam{Focusing on  Witsenhausen's counterexample \cite{Witsenhausen1968} and its variants, 
there are several results in the literature,} which establish lower bounds using information-theoretic techniques (e.g., using the data processing inequality) \cite{bansal1986stochastic, grover2013approximately, Kulkarni2015}, and linear programming-based relaxations \cite{jose2015linear}. However, looking beyond Witsenhausen's counterexample, it is unclear as to how one might extend these techniques to establish computationally tractable lower bounds for the more general family of decentralized control problems considered in this note. More closely related to the approach adopted in this note, there is another stream of literature  that investigates the derivation of computationally tractable lower bounds via information relaxations that increase the amount of information   to which each controller has access to ensure the partial nestedness \cite{asghari2017dynamic, Chu1972, mahajan2010measure, srikant1992asymptotic, yuksel2009stochastic} or quadratic invariance \cite{rotkowitz2012nearest} of the expanded information structure.  }

The specific setting that we consider entails the decentralized output feedback control of a discrete-time, linear time-varying system over a finite horizon, subject to polyhedral  constraints on the state and input trajectories. The system being controlled is partitioned into $N$ dynamically coupled subsystems, each of which has a dedicated local controller. In this setting, the decentralization of information is expressed according to \emph{sparsity constraints} on the information that each local controller has access to. Namely, each local controller is allowed access to the outputs of some subsystems but not others. Naturally, information constraints of this form  may yield information structures that are \emph{nonclassical} in nature, thereby making  the calculation of optimal decentralized control policies  computationally intractable for such systems.
As a result, significant research effort has been directed towards the development of methods for the tractable calculation of constraint admissible, albeit suboptimal, policies. See \cite{bemporad2010decentralized} for an overview.
The aim of this note---which serves to complement these research efforts---is the tractable evaluation of their suboptimality.

\emph{Summary of Results:} \ In this note,  we develop a  tractable approach to the  computation of tight lower bounds on the minimum cost of constrained decentralized control problems with nonclassical information structures. 
The proposed approach is  predicated on two relaxation steps, which together yield  a finite-dimensional convex programming relaxation of the original problem.
The first step entails an information relaxation, which eliminates the so-called signaling incentive between controllers by expanding the set of measurements that each controller has access to. 
\rone{Specifically, we characterize an expansion of the given information structure, which ensures its partial nestedness, while maximizing the optimal value of the resulting decentralized control problem  under the expanded information structure. 
The relaxation is also shown to be \emph{tight}, in the sense that the lower bound induced by the relaxation is achieved for several families of decentralized control problems with nonclassical information.
The relaxed decentralized control problem is then recast as an equivalent convex, infinite-dimensional program using  a nonlinear change of variables akin to the Youla parameterization \cite{Youla1976}.}
Although convex, the resulting optimization problem remains computationally intractable due to its infinite-dimensionality.
As part of the second relaxation step, we obtain a finite-dimensional relaxation of this problem through its partial dualization, and restriction to affine dual control policies.
The resulting problem  is a finite-dimensional  conic program, whose optimal value is guaranteed to be a lower bound on the minimum cost of the original  decentralized control design problem. 
\rone{To the best of our knowledge, such result is the first to offer an efficiently computable (and nontrivial) lower bound on the optimal cost of a decentralized control design problem with multiple subsystems, multiple time periods, and polyhedral constraints on state and input.} If the gap between the cost incurred by an admissible policy and the proposed lower bound is small, then one may conclude that said policy is near-optimal.

\emph{Notation:}  \ Let $\RR$ denote the set of real numbers.
Denote the transpose of a vector $x \in \RR^n$ by $x\trans$.
For any pair of vectors $x = (x_1, .., x_n) \in \RR^n$ and $y = (y_1, .., y_m) \in \RR^m$, we define their concatenation as $(x,y) = (x_1, .., x_n, y_1, .., y_m) \in \RR^{n+m}$. Given a   process $\{x(t)\}$ indexed by  $t = 0, \dots, T-1$, we denote by $x^t = (x(0), x(1), \dots, x(t))$ its history up until and including time $t$. 
We consider block matrices throughout the paper. Given a block matrix $A$ whose  dimension will be clear from the context, we denote by $[A]_{ij}$  its $(i,j)^\text{th}$ block. We denote the trace of a square matrix $A$ by $\trace{A}$.
We denote by $\Kcal_2$ a second-order cone, whose dimension will be clear from the context. \rone{Specifically, given a vector $x \in \RR^n$, $x \in \Kcal_2$ if and only if $x_1 \geq \sqrt{x_2^2 + \cdots + x_n^2}$.} \rsam{Given a  matrix $A$, we let
$A \succeq_{\Kcal_2} 0$  denote its 
columnwise inclusion in $\Kcal_2$. }

\section{Problem Formulation} \label{sec:formulation}

\subsection{System Model}

Consider a discrete-time, linear time-varying system consisting of  $N$ coupled subsystems whose dynamics are described by
\begin{align} \label{eq:x_i}
x_i(t+1) = \sum_{j=1}^N \Big( A_{ij}(t) x_j(t) + B_{ij}(t) u_j(t) \Big) + G_i(t) \xi(t),
\end{align}
 for $i = 1,\dots, N$. The system operates for finite time indexed by $t = 0,\dots, T-1$, and the initial condition is assumed fixed and known.
 We associate with each subsystem $i$  a \emph{local state} $x_i(t) \in \RR^{n_x^i}$ and \emph{local input} $u_i(t) \in \RR^{n_u^i}$. And  we denote by  $\xi(t) \in \RR^{n_\xi}$ the stochastic \emph{system disturbance.}
 We denote by $y_i(t) \in \RR^{n_y^i}$ the \emph{local measured output} of subsystem $i$ at time $t$. It is given by
 \begin{align} \label{eq:y_i}
 y_i(t) = \sum_{j=1}^NC_{ij}(t) x_j(t) + H_i(t) \xi(t),
 \end{align}
for $i = 1,\dots, N$. All system matrices are assumed to be real and of compatible dimension.  \rsam{In the sequel, we will work with a more compact representation of the system Eqs. \eqref{eq:x_i} and \eqref{eq:y_i} given by}
\begin{align*}
x(t+1) & = A(t) x(t) + B(t) u(t) + G(t) \xi(t)\\
y(t) & = C(t) x(t) + H(t) \xi(t).
\end{align*}
Here, we denote by $x(t) := (x_1(t), .., x_N(t)) \in \RR^{n_x}$, $u(t) := (u_1(t), .., u_N(t)) \in \RR^{n_u}$, and $y(t) : = (y_1(t), .., y_N(t)) \in \RR^{n_y}$ the full system state, input, and output at time $t$, respectively. 
\rsam{Their dimensions are given by $n_{x} : = \sum_{i=1}^N n_{x}^i$, $n_{u} : = \sum_{i=1}^N n_{u}^i$, and $n_{y} : = \sum_{i=1}^N n_{y}^i$, respectively.}
We will occasionally refer to  the tuple  $$\Theta := \left\{ A(t),B(t),G(t), C(t), H(t) \right\}_{t=0}^{T-1}$$ as the system parameter when making reference to the underlying system described by Eqs. \eqref{eq:x_i} and \eqref{eq:y_i}.
The system trajectories are related according to
 \begin{align*}
 x = Bu + G\xi \quad \text{and} \quad 
 y  = Cx + H \xi,
 \end{align*}
where $x,u,\xi,$ and $y$ denote the trajectories of the full system state, input, disturbance, and output, respectively.\footnote{\rsam{It is straightforward to construct the block matrices $(A(t),B(t), G(t)$, $C(t), H(t))$  from the  data defining  the system Eqs. \eqref{eq:x_i} and \eqref{eq:y_i}. 
The specification of  the matrices $(B,G,C,H)$ can be found in Appendix \ref{app:matrices}.}} 
We denote them by
\begin{align*}
x & := (x(0), \dots, x(T))  \in \RR^{N_x}, && \rone{ N_x := n_x(T+1),} \\ 
u & := (u(0), \dots, u(T-1))  \in \RR^{N_u}, && \rone{N_u := n_u T,}\\
\xi & : = (1, \xi(0), \dots, \xi(T-1))  \in \RR^{N_\xi}, &&\rone{N_\xi := 1 + n_\xi T,}\\
y & : = (1, y(0), \dots, y(T-1))  \in \RR^{N_y}, && \rone{N_y := 1 + n_y T.}
\end{align*}
Notice that in our specification of the both the  disturbance and output trajectories, $\xi$ and $y$, we have extended each trajectory to include a constant scalar as its initial component. This notational convention will prove useful in simplifying the specification of  affine control policies in the sequel.

\rtwo{
We close this subsection  by stating a structural assumption on the system dynamics. Assumption \ref{ass:prec}, which is assumed to hold throughout the paper,  ensures that each subsystem's local control input can causally  affect its local measured output.
\begin{assumptio} \label{ass:prec}
 For each subsystem $i \in \Vcal$, there exist time periods $0 \leq s < t \leq T-1$ such that the matrix $\left[ C(t) A_{s+1}^t B(s) \right]_{ii}$ is nonzero.
\end{assumptio}
The matrix $\left[C(t) A_{s+1}^t B(s) \right]_{ii}$ refers to 
the $(i, i)^\text{th}$ block of the $N \times N$ block matrix $C(t) A_{s+1}^t B(s)$.  We refer the reader to Appendix \ref{app:matrices} for a definition of the matrix $A_{s+1}^t$.
}

\subsection{Disturbance Model}

We model the disturbance trajectory $\xi$ as a random vector defined according to the probability space $(\RR^{N_\xi}, \Bcal (\RR^{N_\xi}), \PP)$. Here, the Borel $\sigma$-algebra $\Bcal (\RR^{N_\xi})$ denotes the set of all events that are assigned probability according to the measure $\PP$. We denote by $\Lcal^2_n := \Lcal^2 (\RR^{N_\xi}, \Bcal (\RR^{N_\xi}), \PP; \RR^n)$  the space of all $\Bcal (\RR^{N_\xi})$-measurable, square-integrable random vectors taking values in $\RR^n$. 
With a slight abuse of notation, we occasionally use $\xi$ to denote a realization of the random vector $\xi$.  The following 
assumption on the probability distribution of the disturbance trajectory is assumed to hold throughout the paper.
\begin{assumptio}[Elliptically Contoured Disturbance]  \label{ass:noise}
The disturbance trajectory $\xi$ is assumed to have an elliptically contoured distribution.
\rone{That is, there exists a vector $\mu \in \RR^{N_\xi}$, a positive semidefinite matrix $\Sigma \in \RR^{N_\xi \times N_\xi}$, and a scalar function $g$, such that the characteristic function $\varphi_{\xi - \mu}$ of  the random vector $\xi - \mu$ satisfies the functional equation
$\varphi_{\xi - \mu} ( \theta) = g ( \theta\trans \Sigma  \theta)$
for every vector $\theta \in \RR^{N_{\xi}}$. }
\end{assumptio}

The family of elliptically contoured distributions is broad. It includes the multivariate Gaussian distribution, multivariate $t$-distribution, their truncated versions, and uniform distributions on ellipsoids.
If $\xi$ has an elliptically contoured distribution, then  the conditional expectation of $\xi$ given a subvector of $\xi$ is affine in this subvector. And any linear transformation of $\xi$ also follows an elliptically contoured distribution \cite{Cambanis1981}. 
Such properties will play an integral role in the derivation of our main result in Section \ref{sec:dual_LB}.

In order to ensure the well-posedness of the problem to follow, we require that the disturbance trajectory satisfy the following conditions. We assume that the disturbance $\xi$ has support that  is an ellipsoid with a nonempty relative interior in the hyperplane $\{\xi \in \RR^{N_\xi} \, | \, \xi_1 = 1\}$, and is representable by
\begin{align*}
\Xi = \{ \xi \in \RR^{N_\xi} \ | \ \xi_1 = 1 \ \text{and} \  W \xi \succeq_{\Kcal_2} 0 \},
\end{align*}
where $W \in \RR^{N_\xi \times N_\xi}$. This assumption ensures that the corresponding second-order moment matrix  $M := \EE\big[\xi \xi\trans\big]$ is positive definite and finite-valued.

\subsection{System Constraints}

In characterizing the set of feasible input trajectories, we require that the input  and  state trajectories respect the following  linear inequality constraints $\PP$-almost surely,
\begin{align}
\left. \begin{array}{l}
F_x x + F_u u + F_\xi \xi + s = 0 \\
s \geq 0 \\
\end{array} \right\} \ \PP \text{-a.s.}, \label{eq:robust_con2}
\end{align} 
where $F_x \in \RR^{m \times N_x}$, $F_u \in \RR^{m \times N_u}$, and $F_\xi \in \RR^{m \times N_\xi}$.
Here, $s \in \Lcal^2_m$ is a slack variable that is required to be  non-negative $\PP$-almost surely.

\subsection{Decentralized Control Design}

In this paper, we consider information structures that are specified via \emph{sparsity constraints} on the local controllers. More specifically, we describe the pattern according to which information is shared between subsystems with a directed graph $\Gcal_I = (\Vcal,\Ecal_I)$, which we refer to as the \emph{information graph} of the system. Here, the node set $\Vcal = \{1, \dots, N\}$ assigns a distinct node $i$ to each subsystem $i$, and the directed edge set  $\Ecal_I$ determines the pattern of information sharing between subsystems. 
More precisely, we let  $(i,j) \in \Ecal_I$ if and only if for each time $t$, subsystem $j$ has access to subsystem $i$'s local output  $y_i(t)$. We make the following assumption on the structure of the information graph, which ensures that each subsystem $i$ has access to its local output  $y_i(t)$ at each time period $t$.
\begin{assumptio} \label{ass:self}
The directed edge set $\Ecal_I$ is assumed to contain the self-loop $(i,i)$ for each $i \in \Vcal$.
\end{assumptio}
We also assume that each subsystem has \emph{perfect recall}, i.e., each subsystem has access to its entire history of past information at any given time. 
Accordingly, we define the \emph{local information} available to each subsystem $i$ at time $t$ as
\begin{align} \label{eq:pol1}
z_i(t) : = \{ y_j^t \ | \ (j, i) \in \Ecal_I \}.
\end{align}
We restrict the local  input to subsystem $i$ to be of the form
\begin{align} \label{eq:pol2}
u_i(t)  = \gamma_i(z_i(t), t),
\end{align}
where $\gamma_i(\cdot, t)$ is a \rone{measurable function of the local information $z_i (t)$.}
We  define the \emph{local control policy} for subsystem $i$ as $\gamma_i : = (\gamma_i(\cdot,0), \dots, \gamma_i(\cdot,T-1))$. We refer to the  collection of local control policies  $\gamma : = (\gamma_1, \dots, \gamma_N)$ as the \emph{decentralized control policy} and define $\Gamma(\Gcal_I)$ as the set of all decentralized control policies respecting the information structure defined by the information graph $\Gcal_I$.

Of interest is the characterization of  control policies, which solve the following constrained  \emph{decentralized control design problem}:
\begin{equation}
\begin{alignedat}{8}
&\text{minimize}  \ \ &&  \EE \left[ x\trans R_xx + u\trans R_u u\right]\\
& \text{subject to}  \ \ &&  \gamma \in \Gamma(\Gcal_I), \ s \in \Lcal^2_m\\
&&&  \hspace{-.105in} \left. \begin{array}{l}
F_x x + F_u u + F_\xi \xi +s =  0 \\
 x  = Bu + G\xi \\ 
 y  = Cx + H \xi \\
 u = \gamma(y) \\ 
 s \geq 0 \\
\end{array} \right\} \ \PP \text{-a.s.} 
\end{alignedat} \label{opt:decent_full}
\end{equation}
\rone{Here, the \emph{cost matrices}, $R_x \in \RR^{N_x \times N_x}$ and $R_u \in \RR^{N_u \times N_u}$, are both assumed to be symmetric positive semidefinite.}  
We denote the \emph{optimal value} of problem \eqref{opt:decent_full} by $J^* (\Gcal_I) $, \rone{where we have made explicit the dependence of the optimal value of problem \eqref{opt:decent_full} on the underlying information graph $\Gcal_I$.}
In general, the decentralized control design problem \eqref{opt:decent_full} amounts to an infinite-dimensional, nonconvex optimization problem with neither analytical nor computationally efficient solution available at present time \cite{Mahajan2012_survey, Sandell1978, Tsitsiklis1985}.  As a result, considerable effort has been directed towards the development of methods that enable the tractable calculation of feasible control policies \cite{bemporad2010decentralized}. 
Although these methods are known to produce decentralized controllers that perform well empirically,  they are suboptimal in general; and the question as to how far from optimal these policies are remains unanswered. 
The primary objective of this note  is the development of tractable computational methods  to  estimate their suboptimality.

\section{Preliminaries} \label{sec:PN_info}

In what follows, we describe how to equivalently reformulate the decentralized control design problem \eqref{opt:decent_full} as a static team problem \cite{Ho1972} through a nonlinear change of variables akin to the Youla parameterization.  This reformulation is shown to result in a convex program if and only if the underlying  information structure is partially nested.

\subsection{Nonlinear Youla Parameterization} \label{sec:youla}

\rone{Define the  nonlinear Youla parameterization of the decentralized control  policy $\gamma \in \Gamma(\Gcal_I)$ as 
\begin{align}
\phi := \gamma \circ (I - CB \gamma)^{-1}. \label{eq:Youla}
\end{align}
Note that the  map $I - CB \gamma : \RR^{N_y} \to \RR^{N_y}$ is guaranteed to be invertible, as the decentralized control policy $\gamma$ is causal, and the matrix $CB$ is strictly  block lower  triangular. 

The Youla parameter $\phi$ satisfies the following two important properties. First, it is an invertible function of the policy $\gamma$ over $\Gamma(\Gcal_I)$, where its inverse is  given by $\gamma = \phi \circ (I + CB \phi)^{-1}.$  Note that the required  inverse exists, as it is straightforward to verify that $I + CB \phi = (I - CB \gamma)^{-1}$.
Second, given an input trajectory induced by $u =  \gamma(y)$, it holds that 
\begin{align} \label{eq:youlaequiv}
\phi (\eta)=  \gamma (y)
\end{align} 
for every disturbance trajectory $\xi \in \Xi$. Here,  $\eta$  denotes the so-called \emph{purified output trajectory} defined by $\eta := P\xi$, where the matrix  $P \in \RR^{N_y \times N_{\xi}}$ is given by $P := CG + H$.
Note that Eq. \eqref{eq:youlaequiv} follows from the fact that the output trajectory $y$ and purified output trajectory $\eta$ are related according to $$y = CB\gamma(y) + \eta,$$
which in turn implies that  $y = (I - CB \gamma)^{-1}(\eta)$.

Together, these two properties reveal that problem \eqref{opt:decent_full} can be equivalently reformulated as a static team problem by applying the nonlinear change of variables in \eqref{eq:Youla}. This yields the following optimization problem:
\begin{equation}
\begin{alignedat}{8}
&\text{minimize}  \ \ &&  \EE \left[ x\trans R_xx + u\trans R_u u\right]\\
& \text{subject to}  \ \ &&  \phi \in \Phi(\Gcal_I), \ s \in \Lcal^2_m\\
&&&  \hspace{-.105in} \left. \begin{array}{l}
F_x x + F_u u + F_\xi \xi +s =  0 \\
 x  = Bu + G\xi \\ 
 \eta  =P \xi \\
 u =  \phi  (\eta) \\ 
 s \geq 0 \\
\end{array} \right\} \ \PP \text{-a.s.} 
\end{alignedat} \label{opt:decent_full_purified}
\end{equation}
Here,  the set of admissible Youla parameters is given by
$$\Phi (\Gcal_I) := \{   \gamma \circ (I - CB \gamma)^{-1} \ | \  \gamma \in \Gamma (\Gcal_I) \}.$$
The only potential source of  nonconvexity in problem \eqref{opt:decent_full_purified} is in the set of Youla parameters $\Phi (\Gcal_I)$. In particular, problem \eqref{opt:decent_full_purified} is a convex program if and only if the set $\Phi (\Gcal_I)$ is convex. }

\subsection{Convexity under Partially Nested Information Structures}

\rone{In what follows, we show that the static team problem \eqref{opt:decent_full_purified} is a convex program if and only if the information structure is partially nested.}
Before proceeding, we provide a formal definition of partially nested information structures using the notion of precedence, as defined by Ho and Chu in \cite{Ho1974}.
\begin{definitio}[Precedence] \label{def:prec}
Given the information structure defined by $\Gcal_I$, we say subsystem $j$ is a \emph{precedent} to subsystem $i$, denoted by $j \prec i$, if there exist  times  $0 \leq s < t \leq T-1$ and subsystem $k$ satisfying $(k, i) \in \Ecal_I$, such that $\left[ C(t) A_{s+1}^t B(s) \right]_{kj} \neq 0$.
\end{definitio}

Essentially, subsystem $j$ is a \emph{precedent} to subsystem $i$, if the local input to subsystem $j$ can affect the local information available to subsystem $i$ at some point in the future. \rtwo{In particular, it follows from Assumption \ref{ass:prec} that $j$ is a precedent to $i$    if $(j, i) \in \Ecal_I$.} Equipped with the concept of precedence, we now provide the definition of partially nested information structures.

\begin{definitio}[Partially Nested Information]
The information structure defined by $\Gcal_I$ is said to be \emph{partially nested} with respect to the system $\Theta$, if   
$j \prec i$ implies that $z_j (t) \subseteq z_i (t)$
for all times $t=0, \dots, T-1$.
\end{definitio}

We denote by PN$(\Theta)$ the set of information graphs that are partially nested with respect to the the system $\Theta$. 
The information structure defined by $\Gcal_I$ is said to be  \emph{nonclassical} if $\Gcal_I \notin \text{PN} (\Theta)$.
We note that the above definition of partial nestedness is tailored to the setting in which controllers are subject to sparsity constraints on the measured outputs that each controller can access. 
A more general definition of partial nestedness can be found in \cite{Ho1972, Ho1974, Gattami2007}, which applies to the setting in which controllers are subject to both delay and sparsity constraints on information sharing.
\rone{Equipped with this definition, we state the following result, which shows that  the set of Youla parameters $\Phi (\Gcal_I)$ is convex if and only if the information structure is partially nested. We omit the proof of Lemma \ref{lem:PN_necessity}, as it directly follows from existing arguments in \cite[Thm. 1]{rotkowitz2008information} and \cite[Cor. 7]{lessard2011quadratic}. }

\begin{lemm} \label{lem:PN_necessity}
The following statements are equivalent:
\begin{enumerate}[(i)]
\item $ \Phi (\Gcal_I) $ is a convex set,
\item $\Phi (\Gcal_I) = \Gamma (\Gcal_I )$, 
\item $\Gcal_I \in \text{PN} (\Theta)$.
\end{enumerate}
\end{lemm}
\rone{Lemma  \ref{lem:PN_necessity} implies  Ho and Chu's classical result \cite[Thm. 1]{Ho1972} showing that a dynamic team problem with a partially nested information structure can be equivalently reformulated as a static team problem with the same set of admissible policies.}
It follows from \rone{Lemma \ref{lem:PN_necessity}}  that the \rsam{reformulated decentralized control problem} in \eqref{opt:decent_full_purified} is convex if and only if the underlying information structure is partially nested.\footnote{\rone{We note that the convexity result in Lemma \ref{lem:PN_necessity} does not depend on the structure of the cost matrices or the probability distribution of system disturbance. 
There is a related  literature, which identifies structural conditions on the system and cost matrices and the probability distribution of system disturbance, under which the communication of private information from any controller's precedent to said controller does not lead to a reduction in cost. Under these conditions,  the optimal solution of problem \eqref{opt:decent_full_purified} can be computed via the solution of a convex program when the information structure is nonclassical. See \cite{yuksel2009stochastic, yuksel2013stochastic, yuksel2017convex, asghari2017dynamic}  for recent advances.}}

\section{A Convex Information Relaxation} \label{sec:info_relax}

In what follows, we consider systems with nonclassical information structures, and address the question as to how one might convexify  the corresponding decentralized control design problems via information-based relaxations.
Specifically, we \rone{characterize an expansion of the given information graph that guarantees the  partial nestedness of the relaxed information structure, while maximizing the optimal value of the relaxed problem.}
\rone{It is given by the optimal solution to:
\begin{equation}
\begin{alignedat}{8}
& \underset{\Gcal \supseteq \Gcal_I }{\text{maximize}} \quad && \rone{ J^* (\Gcal )} 
 \quad &&\text{subject to} \quad  \Gcal \in \text{PN} (\Theta).
\end{alignedat}\label{opt:min_perturb}
\end{equation}
 Recall that $J^*(\Gcal)$ is the optimal value of the decentralized control design problem \eqref{opt:decent_full} given an information graph $\Gcal$. Also, note that any feasible solution to problem \eqref{opt:min_perturb} is required to both induce a partially nested information structure, and  be a supergraph of $\Gcal_I$.} We  require a few definitions before stating the  solution to problem \eqref{opt:min_perturb}.

\begin{definitio}[Precedence Graph]
We define the \emph{precedence graph} associated with the system $\Theta$ and the information graph $\Gcal_I$ as the directed graph $\Gcal_P (\Theta, \Gcal_I) = (\Vcal, \Ecal_P (\Theta, \Gcal_I) )$ whose directed edge set is defined as $$\Ecal_P (\Theta, \Gcal_I) : = \{(i, j) \; | \; i, j \in \Vcal, \ i \prec j \ \text{with respect to } (\Theta,  \Gcal_I)\}.$$
\end{definitio}
Essentially, the precedence graph provides a directed graphical representation of the precedence relations between all subsystems, as specified in Definition \ref{def:prec}.
\begin{definitio}[Transitive Closure]  The \emph{transitive closure} of a directed graph $\Gcal = (\Vcal, \Ecal)$ is defined as the directed graph $\overline{\Gcal} = (\Vcal, \overline{\Ecal})$, where $(i,j) \in \overline{\Ecal}$ if and only if there exists a directed path in $\Gcal$ from node $i$ to node $j$. 
\end{definitio}
The transitive closure of a directed graph can be efficiently computed using Warshall's algorithm \cite{warshall1962theorem}.
Equipped with these definitions, we state the following result, which provides a `closed-form' solution to problem \eqref{opt:min_perturb}. 
\begin{theore}[Information Relaxation] \label{prop:min_perturbation}
An optimal solution to \eqref{opt:min_perturb} is given by $\overline{\Gcal_P (\Theta, \Gcal_I)}$, the \emph{transitive closure of the precedence graph}.
\end{theore}
Theorem \ref{prop:min_perturbation} implies the following lower bound on the optimal value of the original decentralized control problem \eqref{opt:decent_full}:   
\begin{align} \label{eq:lowerb}
J^* \big( \overline{ \Gcal_P (\Theta, \Gcal_I )} \big) \leq J^*(\Gcal_I).
\end{align} 
Moreover, this  lower bound  can be computed via the solution of a convex infinite-dimensional optimization problem \eqref{opt:decent_full_purified}.
In Theorem \ref{prop:dual_affine_PN}, we  provide a finite-dimensional  relaxation of problem \eqref{opt:decent_full_purified} to enable the tractable approximation of the corresponding lower bound.

\rone{
It is also worth noting that the transitive closure of the precedence graph induces an information structure under which each subsystem is guaranteed to have access to the information of those subsystems whose control input can directly or indirectly affect its information. This implies that the information relaxation $\overline{ \Gcal_P (\Theta, \Gcal_I )}$ yields a partially nested information structure---a result that was originally shown in \cite{Chu1972}.
Theorem \ref{prop:min_perturbation} improves upon this result by establishing the optimality of such a relaxation, in the sense that it is shown to  yield the best lower bound among all partially nested information relaxations.}

\rone{
\begin{remar}[Tightness of the Relaxation]
We also note that the information relaxation in Theorem \ref{prop:min_perturbation} is \emph{tight}. \rsam{That is, $J^* \big( \overline{\Gcal_P (\Theta, \Gcal_I )} \big) =  J^*(\Gcal_I)$ for certain families of  nonclassical control problems.}
In particular, it is known that signaling is \emph{performance irrelevant} if the partially nested information relaxation only introduces additional information that is superfluous in terms of cost reduction---i.e., the additional information does not contribute to an improvement in performance.  
For such problems, one can establish the existence of an optimal policy under the partially nested information relaxation that also respects the original (nonclassical) information structure---implying the tightness of the relaxation.  \rtwo{We refer the reader to \cite{mahajan2010measure}, \cite{yuksel2009stochastic}, \cite[Sec. 3.5]{yuksel2013stochastic} for a rigorous explication of such claims.}
It can also be shown that the lower bound \eqref{eq:lowerb} is achieved  by  nonclassical LQG control problems that  satisfy the so-called    \emph{substitutability condition}. See \cite[Sec. 3]{asghari2017dynamic} for a formal proof of this claim.
\end{remar}
}

In Lemma \ref{lem:PN_property}, we present an alternative characterization of partially nested information structures that will prove useful in the  proof of Theorem \ref{prop:min_perturbation}.
\begin{lemm}  \label{lem:PN_property}
$\Gcal \in \text{PN} (\Theta)$ if and only if $\Gcal =\overline{\Gcal_P (\Theta, \Gcal)}$. 
\end{lemm}
The graph theoretic fixed-point condition in Lemma \ref{lem:PN_property} implies that an information structure is partially nested if and only if the given information graph is equal to the transitive closure of the precedence graph that it induces.
We also note that Lemma \ref{lem:PN_property}  is closely related to the graph theoretic necessary and sufficient condition for quadratic invariance presented in \cite{Swigart2009}, which requires that the information graph be equal to its transitive closure, and be a supergraph of the transitive closure of the so-called plant graph.

\vspace{.1in}

\emph{Proof of Lemma \ref{lem:PN_property}:} \  The proof of the ``if'' direction is straigthforward, and is omitted for brevity. 
We  prove the ``only if'' direction of the statement.  Let $\Gcal = (\Vcal, \Ecal )$.  Assume that $\Gcal \in \text{PN} (\Theta)$. 
\rtwo{It follows from Assumption \ref{ass:prec} that $j \prec i$   if $(j, i) \in \Ecal$.
This implies that $\Gcal \subseteq \Gcal_P (\Theta, \Gcal)$, which in turn implies that $\Gcal \subseteq \overline{\Gcal_P (\Theta, \Gcal)}$.}

To finish the proof, we will show that  $\Gcal \supseteq \overline{\Gcal_P (\Theta, \Gcal)}$.
This amounts to showing that $(j, i) \in \overline{\Ecal_P (\Theta, \Gcal)}$ implies that $(j, i ) \in \Ecal.$
Note that $(j, i) \in \overline{\Ecal_P (\Theta, \Gcal)}$ implies that $j$ is path connected to $i$ in the corresponding precedence graph $\Gcal_P (\Theta, \Gcal)$. That is, there exist $m \geq 1$ distinct vertices $v_1, \dots, v_m \in \Vcal$ that satisfy  $j = v_1 \prec v_2 \prec \dots \prec v_m = i.$
Since $\Gcal \in \text{PN} (\Theta)$, it also holds that
$ z_{v_1} (t) \subseteq z_{v_2} (t) \subseteq \dots \subseteq z_{v_m} (t) $
for  each time $t$.
In particular, it holds that  $z_j (t) \subseteq z_i (t)$ for each time $t$.
This nesting of  information, in combination with Assumption \ref{ass:self}, implies that $(j, i) \in \Ecal$. 
It follows that $\Gcal \supseteq \overline{\Gcal_P (\Theta, \Gcal)},$ which completes the proof. $\hfill \blacksquare$

\vspace{.1in} 

\rone{
We have the following  Corollary to Lemma \ref{lem:PN_property} showing that any graph, which is feasible for problem
\eqref{opt:min_perturb}, must also be a supergraph of the transitive closure of the precedence graph. In other words, this result precludes the existence of feasible information graph relaxations that do not contain $\overline{\Gcal_P (\Theta, \Gcal_I)}$ as a subgraph.
\begin{corollar} \label{cor:feasible_sol}
If  $\Gcal \in \text{PN} (\Theta)$ and $\Gcal  \supseteq \Gcal_I$, then $\Gcal \supseteq \overline{\Gcal_P (\Theta, \Gcal_I)} $.
\end{corollar}

\emph{Proof of Corollary \ref{cor:feasible_sol}:} \  
Lemma \ref{lem:PN_property} implies that $\Gcal = \overline{\Gcal_P (\Theta, \Gcal)} $.
The  result follows, as $\Gcal \supseteq \Gcal_I$ implies that $\overline{\Gcal_P (\Theta, \Gcal)} \supseteq \overline{\Gcal_P (\Theta, \Gcal_I)} $.
$\hfill \blacksquare$}

\vspace{.08in} 

\rone{
\emph{Proof of Theorem \ref{prop:min_perturbation}:} \   Corollary \ref{cor:feasible_sol} implies  that $J^*\big(\overline{\Gcal_P (\Theta, \Gcal_I)} \big) \geq J^*(\Gcal)$ for every graph $\Gcal$ that is feasible for problem \eqref{opt:min_perturb}. Hence, to prove the result, it suffices to show that the graph $\overline{\Gcal_P (\Theta, \Gcal_I)}$ is also feasible for problem \eqref{opt:min_perturb}. 
 We previously showed  in the proof of Lemma \ref{lem:PN_property}  that $\overline{\Gcal_P (\Theta, \Gcal_I)} \supseteq \Gcal_I$.  
We complete the proof by showing that $\overline{\Gcal_P (\Theta, \Gcal_I)} \in \text{PN} (\Theta) $. It is not difficult to show that 
\begin{align} \label{eq:fixedpoint}
\overline{\Gcal_P (\Theta, \Gcal_I)} = \overline{\Gcal_P \left(\Theta, \overline{\Gcal_P (\Theta, \Gcal_I)}\right)}.
\end{align}
This follows from the observation that each precedence relation $ i \prec j$ induced  by $\overline{\Gcal_P (\Theta, \Gcal_I)} $ necessarily corresponds to an edge $ (i,j) \in \overline{\Ecal_P(\Theta, \Gcal_I)}$. 
It follows from \eqref{eq:fixedpoint} and Lemma \ref{lem:PN_property}  that  $\overline{\Gcal_P (\Theta, \Gcal_I)} \in \text{PN}(\Theta)$, which completes the proof. $\hfill \blacksquare$
}


\section{A Dual Approach to Constraint Relaxation} \label{sec:dual_LB}

\rone{The information relaxation developed in Section \ref{sec:info_relax} provides a convex programming relaxation of the original decentralized control design problem \eqref{opt:decent_full}.}
Despite its convexity, the resulting optimization problem remains computationally intractable due to its infinite-dimensionality. 
In what follows, we employ a general technique from robust optimization
\cite{kuhn2011primal,georghiou2015generalized, Hadjiyiannis2011} to obtain a finite-dimensional relaxation of this problem  through its partial dualization, and restriction to affine dual control policies.
The resulting problem  is a finite-dimensional conic program, whose optimal value is guaranteed to be a lower bound on the minimum cost of the original  decentralized control design problem \eqref{opt:decent_full}. 

For the remainder of this section, we assume that the given information structure is partially nested, i.e.,  $\Gcal_I \in \text{PN} (\Theta)$.

\subsection{Restriction to Affine Dual Control Policies}
The derivation of our lower bound centers on a partial Lagrangian relaxation of problem  \eqref{opt:decent_full}. 
We do so by  introducing a \emph{dual control policy} $v \in \Lcal^2_m$, and dualizing the linear equality constraints on the state and input trajectories. This gives rise to the following min-max  problem, which is equivalent to  problem \eqref{opt:decent_full}:
\begin{equation}
\begin{alignedat}{8}
&\text{minimize}  \ \ &&  \sup_{v \in \Lcal^2_m }  \EE  \Big[x\trans R_xx + u\trans R_u u  \\
&&&  \qquad \ + \  v\trans(F_x x + F_u u + F_\xi \xi + s) \Big]  \\
&\text{subject to}  \ \ &&  \gamma \in \Gamma(\Gcal_I), \ \  s \in \Lcal^2_m \\
&&&  \hspace{-.105in} \left. \begin{array}{l}
 x  = Bu + G\xi \\ 
 \eta = P \xi \\
 u  = \gamma(\eta) \\
 s \geq 0 \\
\end{array} \right\} \ \PP \text{-a.s.}
\end{alignedat} \label{opt:dual}
\end{equation}
\rtwo{In presenting the equivalent min-max reformulation of problem \eqref{opt:decent_full}, we have used the fact that problem \eqref{opt:decent_full} is equivalent to problem \eqref{opt:decent_full_purified}; and Lemma \ref{lem:PN_necessity}, which implies that $\Phi (\Gcal_I) = \Gamma (\Gcal_I)$ if $\Gcal_I \in \text{PN} (\Theta)$.}

In order to obtain a tractable relaxation of problem \eqref{opt:dual}, we restrict ourselves to dual control policies that are affine in the disturbance trajectory, i.e., $v = V \xi$ for some $V \in \RR^{m \times N_\xi}$. 
With this restriction, it is possible to derive a closed-form solution for the inner maximization in problem \eqref{opt:dual}.  This yields another minimization problem, whose optimal value stands as a lower bound on that of problem \eqref{opt:dual}. We have the following result, which clarifies this claim.

\begin{propositio} \label{prop:relax}
The optimal value of the following problem is a lower bound on the optimal value of problem \eqref{opt:dual}:
\begin{equation}
\begin{alignedat}{8}
&\text{minimize}  \ \ && \sup_{V \in \RR^{m \times N_\xi} }  \EE  \Big[x\trans R_xx + u\trans R_u u  \\
 &&&  \hspace{.57in} + \,  \xi\trans V\trans(F_x x + F_u u + F_\xi \xi + s) \Big]  \\
&\text{subject to}  \ \ &&  \gamma \in \Gamma(\Gcal_I), \ \  s \in \Lcal^2_m \\
&&&  \hspace{-.105in} \left. \begin{array}{l}
 x  = Bu + G\xi \\ 
 \eta = P \xi \\
 u  = \gamma(\eta) \\
 s \geq 0 \\
\end{array} \right\} \ \PP \text{-a.s.} 
\end{alignedat} \label{opt:dual_aff}
\end{equation}
Moreover, the optimal value of problem \eqref{opt:dual_aff} equals that of the following optimization problem:
\begin{equation}
\begin{alignedat}{8}
&\text{minimize}  \ \ &&  \EE \left[x\trans R_xx + u\trans R_u u \right] \\
& \text{subject to}  \ \ &&  \gamma \in \Gamma(\Gcal_I), \ \  s \in \Lcal^2_m \\
&&& \EE  \left[ (F_x x + F_u u + F_\xi \xi + s) \xi\trans\right] = 0 \\
& && \hspace{-.105in} \left. \begin{array}{l}
 x  = Bu + G\xi \\ 
 \eta = P \xi \\
 u  = \gamma(\eta) \\
 s \geq 0 \\
\end{array} \right\} \ \PP \text{-a.s.} 
\end{alignedat}\label{opt:dual_aff_reform}
\end{equation}
\end{propositio}

\emph{Proof of Proposition \ref{prop:relax}:} \  The fact that the optimal value of problem \eqref{opt:dual_aff} lower bounds that of \eqref{opt:dual} is straightforward, since any dual affine control policy $v = V \xi$ is feasible for the inner maximization problem in \eqref{opt:dual}. To see that the optimal values of problem \eqref{opt:dual_aff} and \eqref{opt:dual_aff_reform} are equal, we note that
\begin{align*}
&\sup_{V \in \RR^{m \times N_\xi} } \EE \Big[  \xi\trans V\trans(F_x x + F_u u + F_\xi \xi + s) \Big]  \\
 &\ \ \qquad = \sup_{V \in \RR^{m \times N_\xi} } \EE \Big[ \trace{V\trans(F_x x + F_u u + F_\xi \xi + s) \xi\trans} \Big] \\
& \ \ \qquad = \begin{cases}
0, & \text{if } \EE \left[ (F_x x + F_u u + F_\xi \xi + s) \xi\trans\right] = 0, \\
+ \infty, & \text{otherwise}.
\end{cases}
\end{align*} $\hfill \blacksquare$

\subsection{Relaxation to a  Finite-dimensional Conic Program}

Problem \eqref{opt:dual_aff_reform} appears to be intractable, as it entails the optimization over an infinite-dimensional function space. In what follows, we show that it admits a relaxation in the form of a finite-dimensional conic program. 
Before proceeding, we present a formal definition of the subspace of causal affine controllers respecting the information structure defined by  $\Gcal_I$.

\begin{definitio} \label{def:aff_control}
Define $S (\Gcal_I) \subseteq \RR^{N_u \times N_y}$ to be the linear subspace of all causal affine controllers respecting the information structure  defined by $\Gcal_I$. 
\end{definitio}
In other words,  for all $K \in S (\Gcal_I)$, the decentralized control policy defined by $\gamma (y) := K y$ satisfies $\gamma \in \Gamma (\Gcal_I)$.
Equipped with this definition, we state the following result, which provides a finite-dimensional relaxation of  problem \eqref{opt:dual_aff_reform} as a conic program.  We note that the proposed conic relaxation is largely inspired by the duality-based relaxation methods originally developed in the context of  centralized control design problems \cite{Hadjiyiannis2011, kuhn2011primal}.  We provide a proof of Theorem \ref{prop:dual_affine_PN} in Appendix \ref{pf:prop:dual_affine_PN}, which extends these techniques to accommodate the added complexity of decentralized information constraints on the controller.

\begin{theore} \label{prop:dual_affine_PN}
Let Assumption \ref{ass:noise} hold. If $\Gcal_I \in \text{PN} (\Theta)$, then the optimal value of the following problem is a lower bound on the optimal value of problem \eqref{opt:decent_full}:
\begin{equation}
\begin{alignedat}{8}
& \text{minimize} \quad & & \trace{P\trans Q\trans R QPM + 2G\trans R_x B QPM + G\trans R_x G M}  \\
& \text{subject to} \quad && Q \in S (\Gcal_I), \ \ Z \in \RR^{m \times N_\xi }\\
&& & (F_u + F_x B) QP + F_x G  + F_\xi + Z = 0, \\
&&& W M Z\trans \succeq_{\Kcal_2} 0,   \\
&&& e_1\trans M Z\trans \geq 0, 
\end{alignedat} \label{eqn:dual_affine_PN}
\end{equation}
where $R = R_u + B\trans R_x B$, and $e_1 = (1,0, \dots, 0)$  is a unit vector in $\RR^{N_\xi}$.
\end{theore}

Let  $J^d (\Gcal_I)$ denote the optimal value of the finite-dimensional conic program \eqref{eqn:dual_affine_PN}.  Theorem \ref{prop:dual_affine_PN} states that  $J^d (\Gcal_I) \leq J^* (\Gcal_I)$  if  $\Gcal_I \in \text{PN} (\Theta)$. 
\rthree{The following  result---an immediate corollary to Theorems \ref{prop:min_perturbation} and  \ref{prop:dual_affine_PN}---provides a  computationally tractable lower bound for  problems with nonclassical information structures.}
\begin{corollar} \label{cor:lowerbound} Let $J^* (\Gcal_I)$ denote the optimal value of the decentralized control design problem \eqref{opt:decent_full}. It follows that
\begin{align*}
 J^d (\overline{\Gcal_P (\Theta, \Gcal_I )} )  \leq J^* (\Gcal_I),
\end{align*}
where $\overline{\Gcal_P (\Theta, \Gcal_I )}$ refers to the transitive closure of the precedence graph associated with problem \eqref{opt:decent_full}.
\end{corollar}

\begin{appendices}

\section{Proof of Theorem \ref{prop:dual_affine_PN}} \label{pf:prop:dual_affine_PN}

The   crux of the proof centers on the introduction of new finite-dimensional decision variables, which enable the removal of the infinite-dimensional decision variables in problem \eqref{opt:dual_aff_reform}. Consider the following result, which we prove in Appendix \ref{pf:lem:finite_dim_var}.

\begin{lemm} \label{lem:finite_dim_var}
Let  Assumption \ref{ass:noise} hold. For each $s \in \Lcal_m^2$, there exists a matrix $Z \in \RR^{m \times N_\xi}$ that satisfies
\begin{align}
ZM = \EE \big[ s \xi\trans \big]. \label{eq:ZM}
\end{align}
For each  $\gamma \in \Gamma(\Gcal_I)$, there exists a matrix $Q \in  S( \Gcal_I)$ that satisfies
\begin{align}
QPM = \EE \big[ u \xi\trans \big], \label{eq:QPM}
\end{align}
where $u = \gamma (\eta)$. 
\end{lemm}
With Lemma \ref{lem:finite_dim_var} in hand, we obtain an equivalent reformulation of problem \eqref{opt:dual_aff_reform} as the following optimization problem---via the  introduction of the finite-dimensional decision variables $Z$ and $Q$ through the constraints \eqref{eq:ZM} and \eqref{eq:QPM}, respectively.

\begin{equation}
\begin{alignedat}{8}
&\text{minimize}  \ \ &&  \EE \left[u\trans R u \right] + \trace{2G\trans R_x B QPM + G\trans R_x G M} \\
& \text{subject to}  \ \ & & \gamma \in \Gamma(\Gcal_I), \ \  s \in \Lcal^2_m, \ \ Q \in S (\Gcal_I), \ \ Z \in \RR^{m \times N_\xi } \\
&&  & (F_u + F_x B) QP M + F_x G M + F_\xi M + Z M = 0 \\
&&& QPM = \EE \big[ u \xi\trans \big] \\
&& & ZM = \EE \big[ s \xi\trans \big] \\
&&& \hspace{-.105in} \left. \begin{array}{l}
 \eta = P \xi \\
 u  = \gamma(\eta) \\
 s \geq 0 \\
\end{array} \right\} \ \PP \text{-a.s.}
\end{alignedat}\label{opt:dual_aff_reform2}
\end{equation}
where $R = R_u + B\trans R_x B$. 

We now introduce two technical Lemmas that permit us to construct a  finite-dimensional relaxation of problem \eqref{opt:dual_aff_reform2}.

\begin{lemm}\label{lem:gamma_opt}
\rthree{Fix the matrix $Q \in S(\Gcal_I)$. It follows that $\gamma (\eta) = Q \eta$ is an optimal solution to the following optimization problem:
\begin{equation*}
\begin{alignedat}{8}
&\text{minimize}  \ \ &&  \EE \left[u\trans R u \right] \\
& \text{subject to}  \ \ & & \gamma \in \Gamma (\Gcal_I) \\
&&& QPM = \EE \big[ u \xi\trans \big] \\
&&& \hspace{-.105in} \left. \begin{array}{l}
 \eta = P \xi \\
 u  = \gamma(\eta) \\
\end{array} \right\} \ \PP \text{-a.s.}
\end{alignedat}
\end{equation*}}
\end{lemm}

We omit the proof of Lemma \ref{lem:gamma_opt}, as it is an immediate corollary of \cite[Lem. 4.5]{Hadjiyiannis2011}.
A direct application of Lemma \ref{lem:gamma_opt} yields  the following equivalent reformulation of problem \eqref{opt:dual_aff_reform2} as:
\begin{equation}
\begin{alignedat}{8}
& \text{minimize} \quad & & \trace{P\trans Q\trans R QPM + 2G\trans R_x B QPM + G\trans R_x G M}  \\
& \text{subject to}  \ \ & &  s \in \Lcal^2_m, \ \ Q \in S (\Gcal_I), \ \ Z \in \RR^{m \times N_\xi } \\
&&  & (F_u + F_x B) QP + F_x G  + F_\xi + Z = 0 \\
&& & ZM = \EE \big[ s \xi\trans \big] \\
&&&  s \geq 0 \quad \PP \text{-a.s.}
\end{alignedat} \label{opt:dual_aff_reform3}
\end{equation}
Note that, in reformulating  problem \eqref{opt:dual_aff_reform2}, we have eliminated the second-order moment matrix $M$ from the equality constraint $(F_u + F_x B) QP M + F_x G M + F_\xi M + Z M = 0$, as $M$ is assumed to be positive definite, and, therefore, invertible.

Lemma \ref{lem:s_relax} provides a conic relaxation of the constraints in problem \eqref{opt:dual_aff_reform3} involving the infinite-dimensional decision variable $s \in \Lcal_m^2$.  We provide a proof of this technical Lemma in Appendix \ref{pf:lem:s_relax}.
\begin{lemm} \label{lem:s_relax}
If $s \in \Lcal_m^2$  and $Z \in \RR^{m \times N_\xi}$ satisfy $ZM = \EE [s \xi\trans]$ and $s \geq 0$  $\PP$-a.s., then
 $W M Z\trans \succeq_{\Kcal_2} 0$ and $e_1\trans M Z\trans \geq 0$.
\end{lemm}

\rone{
We complete the proof with the following string of inequalities and equalities relating the optimal values of the various optimization problems formulated thus far. 
\begin{align*}
\eqref{eqn:dual_affine_PN} \underset{(a)}{\leq} \eqref{opt:dual_aff_reform3}   \underset{(b)}{=}  \eqref{opt:dual_aff_reform2}   \underset{(c)}{=}   \eqref{opt:dual_aff_reform} \underset{(d)}{\leq}   \eqref{opt:dual}  \underset{(e)}{=}    \eqref{opt:decent_full_purified} = \eqref{opt:decent_full} 
\end{align*}
Inequality (a) follows from Lemma \ref{lem:s_relax}, which implies that problem \eqref{eqn:dual_affine_PN} is a relaxation of problem \eqref{opt:dual_aff_reform3}. Equality (b) follows from Lemma \ref{lem:gamma_opt}. Equality (c) follows from Lemma \ref{lem:finite_dim_var}. Inequality (d) follows from Proposition \ref{prop:relax}. Finally,  Equality (e) follows from Lemma \ref{lem:PN_necessity}, as the assumption of a partially nested information structure 
implies equivalence between the optimal values of problems \eqref{opt:dual} and \eqref{opt:decent_full_purified}. The equivalence between \eqref{opt:decent_full_purified} and \eqref{opt:decent_full} is argued in Section \ref{sec:youla}.
}

\section{Proof of Lemma \ref{lem:finite_dim_var}} \label{pf:lem:finite_dim_var}

This proof extends arguments originally developed in \cite[Lem 4.4]{Hadjiyiannis2011} to accommodate the more general setting considered in this note, where the affine controller $Q$ is subject to a decentralized information constraint.

\emph{Proof of the first part:} \ Fix  $s \in \Lcal_m^2$. 
The matrix $M$ is invertible, since it is assumed to be positive definite. Setting $Z =  \EE [ s \xi\trans \big] M^{-1}$ yields the desired result in  \eqref{eq:ZM}.

\vspace{.1in} 

 \rone{
\emph{Proof of the second part:} \ We first introduce the notion of a truncation operator.
Given a nonempty set of  indices $J \subseteq \{1, \dots, N_y\}$ and an arbitrary vector $x \in \RR^{N_y}$, we define $x_J \in \RR^{|J|}$ to be the subvector of $x$, whose entries are given by $x_j$ for all $j \in J$. The entries of $x_J$ are ordered in ascending order of their indices. For example,  if $J = \{1, 3\}$, then $x_J = (x_1, x_3)$.
 We define  the \emph{truncation operator} $\Pi_J: \RR^{N_y} \to \RR^{|J|}$  as the  mapping from a vector $x$ to its subvector $x_J$, i.e., $ x_J = \Pi_J x.$

Now, fix  $\gamma \in \Gamma(\Gcal_I)$, and let $u = \gamma(\eta)$. 
The following  Lemma will prove useful in establishing the existence of a matrix $Q \in S (\Gcal_I)$ satisfying Eq. \eqref{eq:QPM}.  The proof of Lemma \ref{lem:q} is in Appendix \ref{pf:lem:q}.
\begin{lemm} \label{lem:q}
Let Assumption \ref{ass:noise} hold. Let $z \in \Lcal_1^2$ be random variable that is a (possibly nonlinear) function of the random vector $\eta_J = \Pi_J \eta$, where $\{1\} \subseteq J \subseteq \{1, \dots, N_y \}$ is a given index set. Then, there exists another random variable $\widetilde{z} \in \Lcal_1^2$, which is an \emph{affine}\footnote{\rone{As a matter of notational convenience, we have required that $1 \in J$. This enables one to  represent affine functions of $\eta_J$ as linear functions of $\eta_J$, since $\eta_1 = 1$ by construction.}} function of $\eta_J$, and satisfies $  \EE \big[ \widetilde{z} \eta\trans\big] = \EE \big[z \eta\trans\big].$
\end{lemm}
Stated in other words, Lemma \ref{lem:q} asserts the existence of a vector $q \in \RR^{N_y}$ that satisfies
\begin{align}
\EE \big[q\trans \eta \eta\trans\big] = \EE \big[z \eta\trans \big],  \label{eq:q}
\end{align}
where the vector $q$ respects the sparsity pattern encoded by the index set $J$, i.e., $q = \Pi_J\trans \Pi_J q$. 
It follows that one can apply  Lemma \ref{lem:q} to each row of the matrix $\EE [u \eta\trans]$ to establish the existence of a matrix  $Q \in S(\Gcal_I)$ that satisfies 
\begin{align}
\EE\big[Q \eta \eta\trans\big] = \EE \big[u \eta\trans \big]. \label{eq:u_eta}
\end{align}
Consider a matrix $Q \in S(\Gcal_I)$ that satisfies   Eq. \eqref{eq:u_eta}. We complete the proof by showing that this matrix also satisfies Eq. \eqref{eq:QPM}.
First note that the combination of Assumption \ref{ass:noise} and \cite[Thm. 1]{Cambanis1981} implies that the random vector $(\xi, \eta) = (\xi, P \xi)$ has an elliptically contoured distribution. Hence, it follows from \cite[Cor. 5]{Cambanis1981} that the conditional expectation of $\xi$ given $\eta$ is an affine function of $\eta$. The definition of the matrix $P$ ensures that $\eta_1 = 1$. Hence, the conditional expectation can be expressed as
\begin{align}
\EE [ \xi | \eta] = L \eta \quad \PP\text{-a.s.} \label{eq:conditional_exp1}
\end{align}
 for some matrix $L \in \RR^{N_{\xi} \times N_y}$. It holds that
\begin{align*}
\EE \big[ u \xi \trans \big] =  \EE \Big[  \EE \big[ u \xi \trans  \big| \eta \big] \Big] = \EE \big[ u  \eta \trans  \big] L \trans = Q \EE \big[ \eta \eta \trans  \big] L \trans.
\end{align*}
Here, the first equality follows from the law of iterated expectations; the second equality follows from the fact that $u = \gamma (\eta)$ and a direct application of Eq. \eqref{eq:conditional_exp1}; and the third equality follows from Lemma \ref{lem:q}. It also holds that 
\begin{align*}
 \EE \big[\eta \eta\trans\big] L\trans  = \EE \Big[ \eta \EE \big[\xi\trans  \big| \eta \big]\Big] = \EE  \Big[\EE\big[\eta \xi\trans \big| \eta\big] \Big ] = \EE \big[\eta \xi\trans \big]= PM, 
\end{align*}
which completes the proof. 
}

\section{Proof of Lemma \ref{lem:s_relax}} \label{pf:lem:s_relax}

It follows from the symmetry of the matrix $M$ that $M Z\trans = (ZM)\trans = \EE [\xi s\trans]$. It, therefore, holds that 
\begin{align*}
e_1\trans M Z\trans  = e_1\trans \EE [\xi s\trans] = \EE [e_1\trans \xi s\trans] = \EE [s\trans] \geq 0.
\end{align*}
The last equality follows from the fact that $e_1\trans \xi = 1$ $\PP$-almost surely. To show that $W M Z\trans \succeq_{\Kcal_2} 0$, it suffices to show columnwise inclusion in the second-order cone, i.e., 
\begin{align*}
W \EE [s_i \xi ] \succeq_{\Kcal_2} 0, \quad \text{for} \ \ i = 1, \dots, m,
\end{align*}
 where $s_i \in \Lcal_1^2$ is the $i^{\text{th}}$ element of the random vector $s$. By definition, we have that $W \xi \succeq_{\Kcal_2} 0$ for all $\xi \in \Xi$. Also, since  $s_i \geq 0$  almost surely, we have that $W (s_i \xi) \succeq_{\Kcal_2} 0$ almost surely. It follows from the convexity of the second-order cone  that $W \EE [s_i \xi ] \succeq_{\Kcal_2} 0$.

\section{Proof of Lemma \ref{lem:q}} \label{pf:lem:q}

\def \pjmhalf{\Psi}

\rone{
Define the  vector $\subvec \in \RR^{|J|}$ according to
\begin{align}
\subvec\trans  := \EE \big[ z  \eta_J\trans \big] \big( P_J M P_J\trans \big)^{\dagger} , \label{eq:subvec}
\end{align}
where $P_J := \Pi_J P$, and   $ ( \cdot )^{\dagger}$ denotes the Moore-Penrose pseudoinverse of a matrix. 
We first show that the vector $\subvec$ satisfies 
\begin{align}
\EE \big[ z \eta_J\trans \big] = \subvec\trans P_J M P_J\trans . \label{eq:z_eta_J}
\end{align}
Define the matrix $\pjmhalf := P_J M^{1/2}$, where $M^{1/2}$ is the unique square root of the symmetric positive definite matrix $M$. Note that the matrix $M^{1/2}$ is symmetric and positive definite (and hence invertible).
It holds that
\begin{align*}
 \subvec\trans P_J M P_J\trans &= \EE \big[ z  \eta_J\trans \big] \big( P_J M P_J\trans\big)^{\dagger} P_J M P_J\trans \\
& = \EE \big[ z  \xi\trans \big] P_J\trans \big( P_J M P_J\trans\big)^{\dagger} P_J M P_J\trans  \\
& = \EE \big[ z  \xi\trans \big]M^{-1/2} M^{1/2} P_J\trans \big( P_J M P_J\trans\big)^{\dagger} P_J M P_J\trans  \\
&  = \EE \big[ z  \xi\trans \big]  M^{-1/2} \pjmhalf \trans \big( \pjmhalf \pjmhalf \trans \big)^{\dagger} \pjmhalf  \pjmhalf \trans \\
&  = \EE \big[ z  \xi\trans \big]  M^{-1/2} \pjmhalf ^\dagger \pjmhalf  \pjmhalf \trans  \\
& = \EE \big[ z  \xi\trans \big]  M^{-1/2} \pjmhalf \trans = \EE \big[ z  \xi\trans \big]  P_J \trans = \EE \big[ z  \eta_J\trans \big]  
\end{align*}
The second and the last equalities both follow from the fact that $\eta_J = \Pi_J P \xi  = P_J \xi$. The fourth equality follows from the definition of the matrix $\pjmhalf$ and the symmetry of the matrix $M^{1/2}$.
The fifth equality follows from the fact \cite[Prop. 3.2]{barata2012moore} that $\pjmhalf \trans \big( \pjmhalf \pjmhalf \trans \big)^{\dagger} = \pjmhalf^{\dagger}$. The sixth equality follows from the fact \cite[Prop. 3.1]{barata2012moore} that $\pjmhalf ^\dagger \pjmhalf  \pjmhalf \trans = \pjmhalf^{\trans}$.  It follows that the vector $\subvec$ satisfies Eq. \eqref{eq:z_eta_J}.

Now, define the random variable $\widetilde{z} : = r\trans \eta_J$. Clearly, $\widetilde{z}$ is an affine function of $\eta_J$. We complete the proof by showing that $\widetilde{z}$ satisfies
\begin{align}
\EE \big[ \widetilde{z} \eta\trans \big] = \EE \big[ z \eta\trans \big]. \label{eq:subvec_moment}
\end{align}
First recall that the random vector $(\xi, \eta)$ is shown to have an elliptically contoured distribution in Appendix \ref{pf:lem:finite_dim_var}. Hence, it follows from \cite[Cor. 5]{Cambanis1981} that the conditional expectation of $\eta$ given $\eta_J$ is affine in $\eta_J$. The assumption that $1 \in J$ guarantees that the first entry of $\eta_J$ equals 1. Hence,  there exists a matrix $L_J \in \RR^{N_y \times |J|}$, such that
\begin{align}
\EE \left[ \eta \left| \eta_J \right. \right] = L_J \eta_J \quad \PP\text{-a.s.}. \label{eq:conditional_expect}
\end{align}
It holds that
\begin{align*}
\EE \big[ z \eta\trans \big]  &= \EE \Big[ \EE \big[  z \eta \trans \big| \eta_J\big] \Big] = \EE \big[ z \eta_J \trans \big] L_J \trans = \subvec \trans P_J M P_J \trans L_J \trans 
\end{align*}
Here, the first equality follows from the law of iterated expectations. The second equality follows from a combination of Eq. \eqref{eq:conditional_expect} and the assumption that $z$ is a function of $\eta_J$. The third equality follows from Eq.  \eqref{eq:z_eta_J}. It also holds that 
\begin{align*}
\subvec \trans P_J M P_J \trans L_J \trans & = \subvec \trans \EE \big[ \eta_J \eta_J \trans \big] L_J \trans = \subvec \trans \EE \Big[ \EE \big[ \eta_J \eta \trans \big| \eta_J \big]\Big] \\
& = \subvec \trans \EE \big[ \eta_J \eta \trans \big] = \EE \big[ \widetilde{z} \eta\trans \big],
\end{align*}
which completes the proof.
}

\section{Matrix Definitions} \label{app:matrices}

The block matrices  $(B,G,C,H)$ are given by:
\begin{align*}
 B & := \bmat{0 \\ A_1^1 B(0) &0 \\ A_1^2 B(0) & A_2^2 B(1) & 0 \\  \vdots &  &  & \ddots \\ \vdots &&&&0 \\ A_1^T B(0) & A_2^T B(1)  & \cdots & \cdots &A_T^T B(T-1)}\\
 G & := \bmat{A_0^0 x (0)   \\ A_0^1 x (0) & A_1^1 G(0)  \\ A_0^2 x (0) & A_1^2 G(0) & A_2^2 G(1)   \\ \vdots &  \vdots  &  & \ddots \\  A_0^T x (0) & A_1^T G(0) & A_2^T G(1)  &   \cdots &A_T^T G(T-1)} \\
 C & := \bmat{0 \\
C(0) & 0 \\
& \ddots & \ddots \\
&& C(T-1) & 0} \\
H &:= \text{diag}(1, H(0), \dots , H(T-1)),
\end{align*}
where $A_s^t := \prod_{r=s}^{t-1} A(r)$ for $s < t$, and $A_t^t = I $.

\end{appendices}

\bibliographystyle{IEEEtran}
\bibliography{decentralized_bib}{\markboth{References}{References}}

\end{document}